\documentclass[11pt]{amsart}
\usepackage[T1]{fontenc}
\usepackage[utf8]{inputenc}
\usepackage{lmodern}
\usepackage{amsmath,amssymb,mathtools}
\usepackage{booktabs}
\usepackage{microtype}
\usepackage{tikz}
\usetikzlibrary{trees}
\usepackage[colorlinks=true,linkcolor=blue!55!black,citecolor=blue!55!black,
            urlcolor=blue!55!black]{hyperref}

\usetikzlibrary{trees,backgrounds}

\newtheorem{theorem}{Theorem}[section]
\newtheorem{proposition}[theorem]{Proposition}
\newtheorem{corollary}[theorem]{Corollary}
\newtheorem{lemma}[theorem]{Lemma}

\theoremstyle{definition}

\usepackage{cleveref}

\theoremstyle{remark}
\newtheorem{remark}[theorem]{Remark}

\newcommand{\SF}{\mathcal{SF}}
\newcommand{\sfi}{\operatorname{sf}}
\newcommand{\sodd}{\operatorname{sodd}}
\newcommand{\odd}{\operatorname{odd}}

\newcommand{\depth}{\operatorname{depth}}

\newcommand{\Tsf}{\mathcal{T}_{\mathrm{SF}}}

\calclayout

\title[The Semi--Fibonacci Partition Tree]{The Semi--Fibonacci Partition Tree}

\author{Andr\'es R. Moreno}
\address{Departamento de Matem\'aticas, Universidad del Cauca, Popay\'an, Colombia}
\email{armoreno@unicauca.edu.co}

\author{Jos\'e L. Ram\'irez}
\address{Departamento de Matem\'aticas, Universidad Nacional de Colombia, Bogot\'a, Colombia}
\email{jlramirezr@unal.edu.co}


\begin{document}

\begin{abstract}
In this work, we introduce a binary tree in which every semi-Fibonacci partition appears exactly once. The left--right address of a vertex gives an explicit formula for the parts of the corresponding partition and encodes statistics such as the largest $2$-adic
valuation among its parts. We construct a sign-reversing involution
that refines Andrews's parity theorem by determining the difference
between the numbers of partitions of a given integer at even and odd depths. Inspired by properties of the Calkin--Wilf tree, we
study several parameters of our tree and use them to obtain new
results on semi-Fibonacci partitions. We also derive generating
functions for partition statistics and determine the average and
extremal partition sizes at each level.
\end{abstract}

\maketitle

\section{Introduction}

A \emph{partition} of a positive integer $n$ is a weakly decreasing sequence $\lambda=(\lambda_1,\lambda_2,\ldots,\lambda_k)$ of positive integers whose sum is $n$. The numbers $\lambda_1,\lambda_2,\ldots,\lambda_k$ are called its \emph{parts}.
We write $|\lambda|$ for the sum of its parts and $\ell(\lambda)$ for its number of parts.

\begin{figure}[ht!]
\centering
\resizebox{\textwidth}{!}{%
\begin{tikzpicture}[
  level distance=1.2cm,
  every node/.style={font=\small, inner sep=1pt},
  edge from parent/.style={draw},
  level 1/.style={sibling distance=9.6cm},
  level 2/.style={sibling distance=4.8cm},
  level 3/.style={sibling distance=2.4cm},
  level 4/.style={sibling distance=1.2cm}
]
\node {$(1)$}
  child {node {$(1^2)$}
    child {node {$(1^3)$}
      child {node {$(1^4)$}
        child {node {$(1^5)$}}
        child {node {$(2^4)$}}
      }
      child {node {$(2^3)$}
        child {node {$(2^3,1)$}}
        child {node {$(3^3)$}}
      }
    }
    child {node {$(2^2)$}
      child {node {$(2^2,1)$}
        child {node {$(2^2,1^2)$}}
        child {node {$(3^2,2)$}}
      }
      child {node {$(3^2)$}
        child {node {$(3^2,1)$}}
        child {node {$(4^2)$}}
      }
    }
  }
  child {node {$(2)$}
    child {node {$(2,1)$}
      child {node {$(2,1^2)$}
        child {node {$(2,1^3)$}}
        child {node {$(3,2^2)$}}
      }
      child {node {$(3,2)$}
        child {node {$(3,2,1)$}}
        child {node {$(4,3)$}}
      }
    }
    child {node {$(3)$}
      child {node {$(3,1)$}
        child {node {$(3,1^2)$}}
        child {node {$(4,2)$}}
      }
      child {node {$(4)$}
        child {node {$(4,1)$}}
        child {node {$(5)$}}
      }
    }
  };
\end{tikzpicture}%
}
\caption{The first five levels of the binary tree of integer
partitions. Exponents indicate multiplicities of parts.}
\label{fig:partitiontree}
\end{figure}
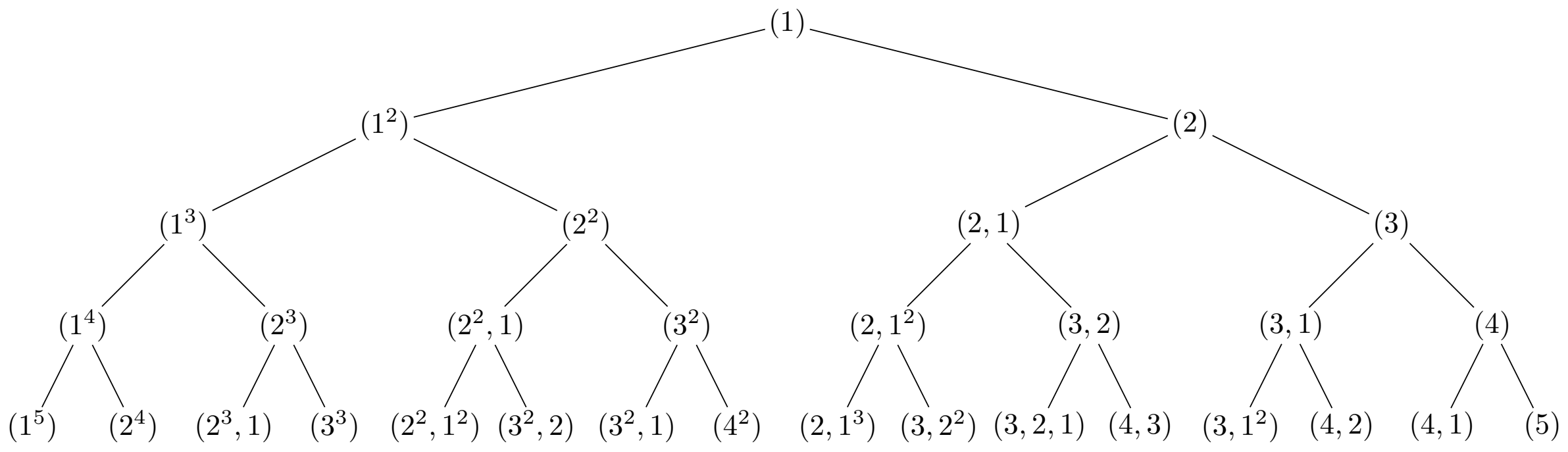

The nonempty integer partitions can be arranged in a binary tree
in which each partition occurs exactly once. 
Place $(1)$ at the root. Given a partition
$\lambda=(\lambda_1,\ldots,\lambda_\ell)$, obtain its left child by
adjoining a part $1$, and its right child by adding $1$ to every part. Thus the two children are $(\lambda_1,\ldots,\lambda_\ell,1)$ and $(\lambda_1+1,\ldots,\lambda_\ell+1)$, respectively. Every partition other than $(1)$ has a unique parent, obtained by
removing a part $1$ if one is present and otherwise subtracting $1$
from every part. Iterating this parent map eventually reaches the root, so every nonempty partition appears exactly once. The first five levels of this partition tree are shown in Figure~\ref{fig:partitiontree}. This tree is closely related to the usual binary encoding of partitions; see \cite[Exercise 7.59]{Stanley2}. Hemmer and Westrem~\cite{HemmerWestrem} used its correspondence with the Calkin--Wilf tree to study palindrome partitions and symmetries of Young diagrams. 
 
We construct an analogous binary tree for \emph{semi-Fibonacci
partitions}, which are defined by the following even--odd recursion. Set $\SF(1)=\{(1)\}$ and $\SF(2)=\{(2)\}$.  For $n\geq1$, the partitions in $\SF(2n)$ are obtained by doubling every part of a partition in $\SF(n)$. The partitions in $\SF(2n+1)$ form the disjoint union of two families: those obtained by adjoining a part $1$ to a partition in $\SF(2n)$, and those obtained by adding $2$ to the unique odd part of a partition in $\SF(2n-1)$. Every partition in $\SF(2n-1)$ has a unique odd part, so the second
operation is well defined. For example,
\begin{align*}
 \SF(3)&=\{(3),(2,1)\}, \qquad  \SF(4)=\{(4)\}, \qquad \SF(5)=\{(5),(4,1),(3,2)\},\\ 
 \SF(6)&=\{(6),(4,2)\}, \qquad \SF(7)=\{(7),(6,1),(5,2),(4,3),(4,2,1)\}.
\end{align*}
Let $\sfi(n)=|\SF(n)|$ for $n\geq1$, and set $\sfi(0)=1$.
The numbers $\sfi(n)$ are called the \emph{semi-Fibonacci numbers}
and satisfy
\[
\sfi(2n)=\sfi(n),\qquad
\sfi(2n+1)=\sfi(2n)+\sfi(2n-1)\qquad(n\geq1).
\]
The first few values of $\{\sfi(n)\}_{n\geq0}$ are
\[
  1,\ 1,\ 1,\ 2,\ 1,\ 3,\ 2,\ 5,\ 1,\ 6,\ 3,\ 9,\ 2,\ 11,\ 5,\ 16,\ldots.
\]
This is sequence A030067 in the OEIS \cite{OEIS}.

Beck~\cite{Beck} considered the subsum polynomial
$P_\lambda(x)=\prod_i(1+x^{\lambda_i})$ associated with a partition
$\lambda$. The coefficient of $x^j$ counts the ways to select
parts of $\lambda$ with sum $j$. Beck proved that $P_\lambda(x)$
has no repeated zeros if and only if $\lambda$ is a semi-Fibonacci
partition. Andrews~\cite{Andrews} proved that, for $n>0$,
$\mathrm{sf}(n)$ is even if and only if $3\mid n$. For further
bijections, generalizations, and congruence properties, see
\cite{Alanazi,Alegri,BallantineBeck,MorenoRamirez}.

Our construction gives an explicit formula for the parts of a partition in terms of its left--right address and translates some partition statistics into statistics of binary words. We number the vertices level by level from left to right, with the root numbered $1$, and let $h_j$ be the size of the partition
at vertex $j$. The equality $h_j=j$ holds if and only if $j=2^a$ or $j=3\cdot2^a$ for some integer $a\ge0$. At every other index,
$h_j/j\le5/7$, and this bound is sharp.

We introduce a sign-reversing involution that proves Andrews's parity theorem and gives a signed refinement by depth. Using the semi-Fibonacci polynomials introduced in our earlier work~\cite{MorenoRamirez}, we determine the first and last levels at which partitions of a fixed integer occur and show
that the gaps between occupied levels can be arbitrarily long.

We also prove that the minimum number of right steps among semi-Fibonacci partitions of $n$ is $s_2(n)-1$, where $s_2(n)$ denotes the binary digit sum of $n$. The partitions attaining this
minimum are obtained by merging disjoint pairs of consecutive powers of two in the binary expansion of $n$, and their
number is a product of Fibonacci numbers. Finally, we derive generating functions for partition statistics and formulas for the average and extremal partition sizes at each depth.

\section{From the recursion to the tree}

We begin with a characterization of semi-Fibonacci partitions in terms of the $2$-adic valuations of their parts. This characterization
gives recursive rules for a binary tree in which every semi-Fibonacci partition occurs exactly once.

For a positive integer $n$, its $2$-adic valuation $\nu_2(n)$ is the largest integer $k\ge0$ such that $2^k\mid n$.

\begin{lemma}[{cf. \cite[Theorem~3.2]{Alanazi}}]\label{lem:valuations}
A partition of $n$ belongs to $\mathrm{SF}(n)$ if and only if its
parts have pairwise distinct $2$-adic valuations.  \end{lemma}

\begin{proof}
We use induction on $n$. The case $n=1$ is immediate. For the forward implication, each operation in the recursive definition preserves distinct valuations. Doubling all parts
increases their valuations by $1$. Adjoining a part $1$ to a partition with only even parts introduces the new valuation $0$, while adding $2$ to the unique odd part leaves its valuation
equal to $0$.

Conversely, let $\lambda$ be a partition of $n$ whose parts have pairwise distinct valuations. Then $\lambda$ has at most one odd
part. If $n$ is even, all parts are even, and dividing them by $2$ gives a smaller partition with pairwise distinct valuations.
If $n$ is odd, let $u$ be the unique odd part. Remove this part if $u=1$, and replace it by $u-2$ otherwise. Again, the resulting
partition has smaller size and pairwise distinct valuations. In either case, the induction hypothesis and the recursive definition give $\lambda\in\mathrm{SF}(n)$.
\end{proof}

For example, $(6,4)\in\SF(10)$, since its parts have $2$-adic valuations $1$ and $2$. On the other hand, $(5,3,2)\notin\SF(10)$, since the parts $5$ and $3$ have the same $2$-adic valuation $0$.

For a partition $\lambda$,  we write $2\lambda$ for the partition obtained by doubling every part
of $\lambda$. We define the \emph{left child} of a semi-Fibonacci partition $\lambda$ by $L(\lambda)=2\lambda$, and its \emph{right child} by
\[
R(\lambda)=
\begin{cases}
\lambda\cup\{1\},&\text{if } |\lambda|\text{ is even},\\
\lambda\text{ with its unique odd part $u$ replaced by $u+2$},
& \text{if } |\lambda|\text{ is odd}.
\end{cases}
\]

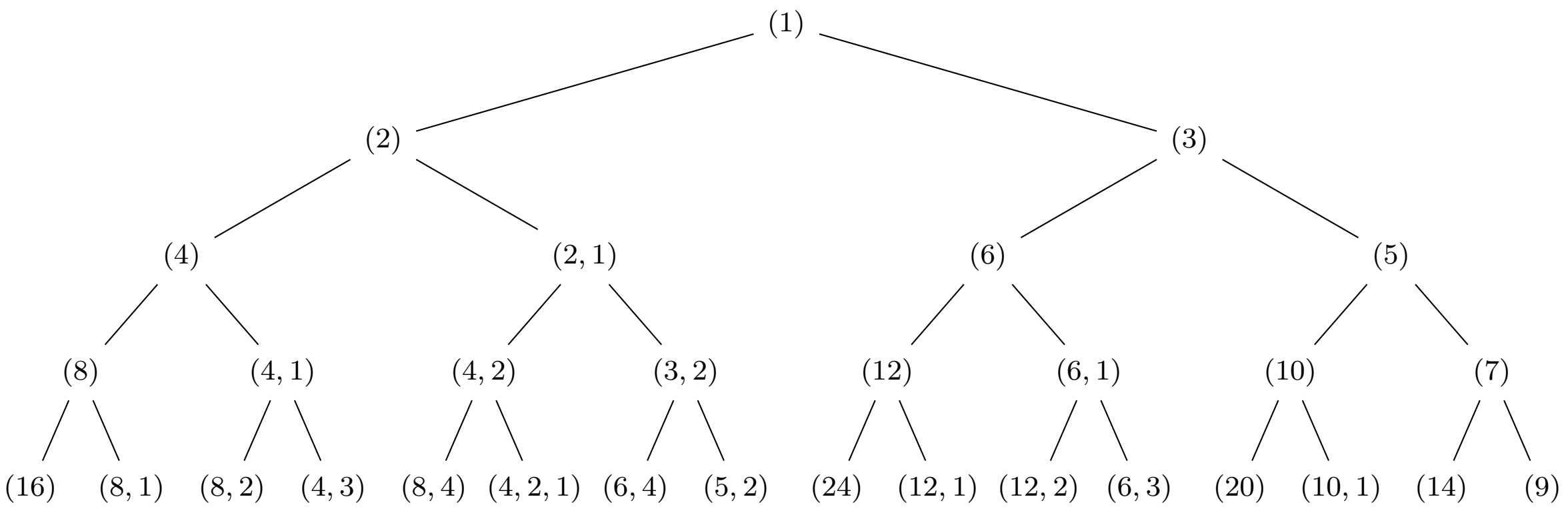
\begin{figure}[ht!]
\centering
\begin{tikzpicture}[
  level distance=1.15cm,
  every node/.style={font=\scriptsize},
  edge from parent/.style={draw},
  level 1/.style={sibling distance=8cm},
  level 2/.style={sibling distance=4cm},
  level 3/.style={sibling distance=2cm},
  level 4/.style={sibling distance=1cm}
]
\node {$(1)$}
  child {node {$(2)$}
    child {node {$(4)$}
      child {node {$(8)$}
        child {node {$(16)$}}
        child {node {$(8,1)$}}
      }
      child {node {$(4,1)$}
        child {node {$(8,2)$}}
        child {node {$(4,3)$}}
      }
    }
    child {node {$(2,1)$}
      child {node {$(4,2)$}
        child {node {$(8,4)$}}
        child {node {$(4,2,1)$}}
      }
      child {node {$(3,2)$}
        child {node {$(6,4)$}}
        child {node {$(5,2)$}}
      }
    }
  }
  child {node {$(3)$}
    child {node {$(6)$}
      child {node {$(12)$}
        child {node {$(24)$}}
        child {node {$(12,1)$}}
      }
      child {node {$(6,1)$}
        child {node {$(12,2)$}}
        child {node {$(6,3)$}}
      }
    }
    child {node {$(5)$}
      child {node {$(10)$}
        child {node {$(20)$}}
        child {node {$(10,1)$}}
      }
      child {node {$(7)$}
        child {node {$(14)$}}
        child {node {$(9)$}}
      }
    }
  };
\end{tikzpicture}
\caption{The first five levels of the semi-Fibonacci partition tree.}
\label{fig:SFtree}
\end{figure}

The \emph{semi-Fibonacci partition tree}, denoted by $\Tsf$,
is the rooted ordered binary tree with root $(1)$. If a vertex is labeled by a partition $\lambda$, then its left child is labeled by $L(\lambda)$
and its right child by $R(\lambda)$. We call an edge from a vertex $\lambda$ to $L(\lambda)$ an \emph{$L$-edge}, and an edge from $\lambda$ to $R(\lambda)$ an \emph{$R$-edge}. The first five levels of $\Tsf$ are shown in Figure~\ref{fig:SFtree}.  

By the same argument used in Lemma~\ref{lem:valuations}, every
semi-Fibonacci partition other than $(1)$ has a unique parent of
smaller size. This gives the following result.

\begin{theorem}\label{thm:partition-tree}
Every semi-Fibonacci partition occurs exactly once in $\Tsf$. Thus, if
$V(\Tsf)$ denotes the vertex set of $\Tsf$, then
\[
V(\Tsf)=\bigcup_{n\geq1}\SF(n).
\]
\end{theorem}

\subsection{Reading a partition from its address}

Every finite word $w$ over $\{L,R\}$ has a unique expression \begin{equation}\label{eq:block-decomposition}
w=R^{r_0}LR^{r_1}L\cdots LR^{r_k},
\qquad r_0,r_1,\ldots,r_k\geq0,
\end{equation}
where $k$ is the number of $L$'s in $w$. The exponents may be zero, so consecutive $L$'s are allowed.

Let $\lambda(w)$ denote the semi-Fibonacci partition reached from the root $(1)$ by following the sequence of $L$- and $R$-steps encoded by $w$.

\begin{proposition}\label{prop:address-formula}
For $w$ as in \eqref{eq:block-decomposition}, the parts of
$\lambda(w)$ are
\begin{equation}\label{eq:partition-from-word}
\bigl\{(2r_0+1)2^k\bigr\}
\cup
\bigl\{(2r_i-1)2^{k-i}:1\leq i\leq k,\ r_i>0\bigr\}.
\end{equation}
In particular,
\[
|\lambda(w)|
=(2r_0+1)2^k
+
\sum_{\substack{1\leq i\leq k\\r_i>0}}
(2r_i-1)2^{k-i}.
\]
\end{proposition}

\begin{proof}
Starting from the root  $(1)$, after following the initial block $R^{r_0}$, the unique odd part is $2r_0+1$. Each of the $k$ occurrences of $L$ doubles this part, so it eventually becomes $(2r_0+1)2^k$. After the $i$th occurrence of $L$, all parts are even. If $r_i>0$, the first $R$ adjoins a part $1$, and 
the remaining $r_i-1$ $R$-steps increase this odd part by $2$ each time.  Thus the block $R^{r_i}$ produces the part $2r_i-1$. The remaining $k-i$ $L$-steps double it to  $(2r_i-1)2^{k-i}$.  If $r_i=0$, no new part is created. This gives
\eqref{eq:partition-from-word}, and the formula for $|\lambda(w)|$
follows by summing the parts.
\end{proof}

For example, for the word $w=LRRL=R^0LR^2LR^0$, we have $(r_0,r_1,r_2)=(0,2,0)$. The corresponding parts are $2^2=4$ and $(2\cdot2-1)2^1=6$, so $\lambda(LRRL)=(6,4)$ and $|\lambda(LRRL)|=10$. Conversely, the address of a semi-Fibonacci partition is recovered by
following its unique parent chain back to $(1)$ and reversing the resulting sequence of $L$- and $R$-steps.

\begin{remark}
The \emph{Calkin--Wilf tree}~\cite{CalkinWilf} has root $1/1$ and child maps
\[
L\left(\frac ab\right)=\frac{a}{a+b},
\qquad
R\left(\frac ab\right)=\frac{a+b}{b}.
\]
Both trees are indexed by finite words over $\{L,R\}$. Matching equal addresses associates a semi-Fibonacci partition with each positive rational. Proposition~\ref{prop:address-formula} and the continued-fraction description of Calkin--Wilf
addresses~\cite{Callan,McConvilleProppSagan} make this
bijection explicit.  For a positive rational $r$, write its  finite regular continued fraction with an odd number of partial quotients as $r=[a_0;a_1,\ldots,a_{2k}]$,
and set
\[
E_0=0,\qquad
E_j=a_1+a_3+\cdots+a_{2j-1}\quad(j\geq1).
\]
The corresponding semi-Fibonacci partition has parts
\[
\bigl\{(2a_{2j}-1)2^{E_j}:
0\leq j\leq k,\ a_{2j}>0\bigr\}.
\]
For example, $\frac57=[0;1,2,1,1]$ corresponds to $(6,4)$.
\end{remark}

For a word $w$ over $\{L,R\}$, let $|w|_L$ and $|w|_R$ denote the numbers of occurrences of $L$ and $R$, respectively, and let $|w|_{LR}$ count the occurrences of the factor $LR$.  We also use two partition statistics introduced in our earlier work~\cite{MorenoRamirez}. For a positive integer $u$, let $\odd(u)$ be its odd part, so that $u=2^{\nu_2(u)}\odd(u)$. For a semi-Fibonacci partition $\lambda=(\lambda_1,\ldots,\lambda_t)$, define
\[
\nu_{\max}(\lambda)=\max_{1\leq i\leq t}\nu_2(\lambda_i) \qquad \text{and} 
\qquad
\sodd(\lambda)=
\sum_{i=1}^{t}
\left\lfloor\frac{\odd(\lambda_i)}{2}\right\rfloor.
\]
The \emph{depth} of $\lambda$, denoted by  $\depth(\lambda)$, is the distance from $\lambda$ to the root $(1)$ in the semi-Fibonacci tree. Thus $\depth(\lambda(w))=|w|$, where $|w|$ is the length of $w$. The next corollary  expresses these partition statistics in terms of the address word

\begin{corollary}\label{cor:word-statistics}
For every finite word $w$ over $\{L,R\}$,
\[
\begin{aligned}
\nu_{\max}(\lambda(w))&=|w|_L,\\
\ell(\lambda(w))&=1+|w|_{LR},\\
\sodd(\lambda(w))&=|w|_R-|w|_{LR}.
\end{aligned}
\]
Consequently, every semi-Fibonacci partition $\lambda$ satisfies
\[
\depth(\lambda)
=
\nu_{\max}(\lambda)+\ell(\lambda)+\sodd(\lambda)-1.
\]
\end{corollary}

\begin{proof}
In \eqref{eq:partition-from-word}, the largest $2$-adic valuation is $k=|w|_L$. In addition to the part arising from $r_0$, there is a part for each positive $r_i$, $i\geq1$. Thus the positive $r_i$ with $i\geq1$ are in bijection with the
occurrences of $LR$ in $w$. Finally,
\[
\sodd(\lambda(w))
=
r_0+\sum_{\substack{1\leq i\leq k\\r_i>0}}(r_i-1)
=
|w|_R-|w|_{LR}.
\]
The identity for the depth now follows from
$|w|=|w|_L+|w|_R=\depth(\lambda(w))$.
\end{proof}

\begin{remark}
For $d\geq1$, choose a vertex $\lambda$ uniformly at random
from the vertices at depth $d$.  Then
\[
\mathbb E\bigl(\nu_{\max}(\lambda)\bigr)=\frac d2,
\qquad
\mathbb E\bigl(\ell(\lambda)\bigr)=\frac{d+3}{4},
\qquad
\mathbb E\bigl(\sodd(\lambda)\bigr)=\frac{d+1}{4}.
\]
Indeed, the address $w$ of $\lambda$ is uniformly distributed
over $\{L,R\}^d$, so
\[
\mathbb E(|w|_L)=\mathbb E(|w|_R)=\frac d2,
\qquad
\mathbb E(|w|_{LR})=\frac{d-1}{4}.
\]
The formulas follow from Corollary~\ref{cor:word-statistics}.
\end{remark}

\section{From the tree to the weight sequence}

Number the vertices of $\Tsf$ level by level from left to right, starting with $1$ at the root. The left and right children of vertex $j$ then have indices $2j$ and $2j+1$, respectively. The \emph{weight} $h_j$ of vertex $j$ is the size of its partition.

Since every semi-Fibonacci partition appears exactly once,
\[
\sfi(n)=\#\{j\geq1:h_j=n\}\qquad(n\geq1).
\]
Thus the semi-Fibonacci numbers count the occurrences of
each positive integer in the weight sequence.

The tree construction gives $h_1=1$ and, for $j\geq1$,
\[
h_{2j}=2h_j,\qquad
h_{2j+1}=
\begin{cases}
h_j+1,&\text{if $h_j$ is even},\\
h_j+2,&\text{if $h_j$ is odd}.
\end{cases}
\]
Left children have even index and even weight, while right children have odd index and odd weight. Hence
$h_j\equiv j\pmod 2$ for every $j\geq1$, and we can rewrite
the second recurrence as
\[
h_{2j+1}=h_j+1+(j\bmod 2).
\]
The first few terms, which can also be read from Figure~\ref{fig:SFtree}, are
\[
1,\ 2,\ 3,\ 4,\ 3,\ 6,\ 5,\ 8,\ 5,\ 6,\ 5,\ 12,\ 7,\ 10,\ 7,\ldots.
\]
A weight can occur at different depths. For example, $h_3=h_5=3$, corresponding to the partitions $(3)$ and $(2,1)$ at depths $1$ and $2$, respectively.

The binary expansion of $j$ records the path to vertex $j$: passing to a left or right child appends a $0$ or $1$, respectively. Deleting the leading $1$ and replacing $0$ and
$1$ by $L$ and $R$ therefore gives its address. Proposition~\ref{prop:address-formula} now gives the following
formula.

\begin{corollary}\label{cor:binary-index}
Let $j$ have binary expansion
\[
j=\bigl(1\,1^{r_0}0\,1^{r_1}0\cdots0\,1^{r_k}\bigr)_2,
\qquad r_0,\ldots,r_k\geq0.
\]
Then
\[
h_j
=
(2r_0+1)2^k
+
\sum_{\substack{1\leq i\leq k\\r_i>0}}
(2r_i-1)2^{k-i}.
\]
\end{corollary}

For example, $22=(10110)_2$ has address $LRRL$, and the formula
gives $h_{22}=10$, corresponding to the partition $(6,4)$.

\begin{figure}[ht!]
    \centering
    \includegraphics[scale=0.6]{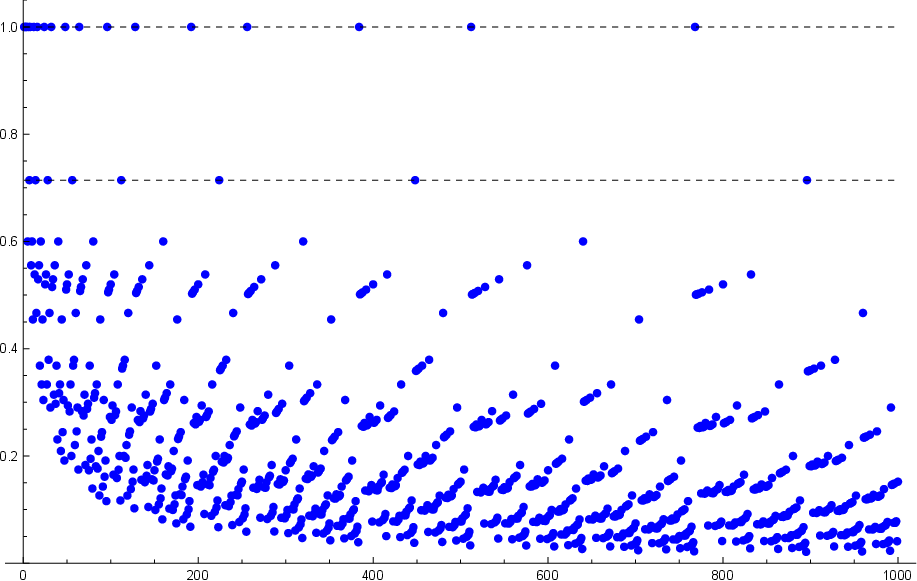}
\caption{The normalized weights $h_j/j$ for $1\leq j\leq1000$.
The dashed lines mark the heights $1$ and $5/7$.}
\label{fig:normalized-weights}
\end{figure}

We next compare the weight $h_j$ with the index $j$. Figure~\ref{fig:normalized-weights} plots the ratios $h_j/j$,
which appear to be either $1$ or at most $5/7$.

We now explain these two features. Define the \emph{defect} $\delta_j=j-h_j$. The recurrence for $h_j$ gives
\[\delta_{2j}=2\delta_j,
\qquad
\delta_{2j+1}
=
\delta_j+2\left\lfloor\frac j2\right\rfloor.
\]
Since $\delta_1=0$, induction gives $\delta_j\geq0$ and hence $h_j\leq j$ for every $j\geq1$.

The next theorem characterizes the indices for which equality holds.

\begin{theorem}\label{thm:mass-fixed-points}
For $j\geq1$, the following conditions are equivalent:
\begin{enumerate}
\item $h_j=j$.
\item $j=2^a$ or $j=3\cdot2^a$ for some integer $a\geq0$.
\item The address of vertex  $j$ is $L^a$ or $RL^a$ for some integer  $a\geq0$.
\item The semi-Fibonacci partition at vertex $j$ is the one-part partition $(j)$.
\end{enumerate}
\end{theorem}

\begin{proof}
First consider an odd index $j=2m+1>1$. Then
\[
\delta_j
=
\delta_m+2\left\lfloor\frac m2\right\rfloor .
\]
Both terms on the right are nonnegative. Since $m\geq1$, the second vanishes only when $m=1$. Together with $\delta_1=0$,
this shows that the only odd fixed points are $1$ and $3$. For even indices, $\delta_{2m}=0$ if and only if $\delta_m=0$. Repeatedly dividing by $2$ therefore gives
\[
\delta_j=0
\quad\Longleftrightarrow\quad
j=2^a
\quad\text{or}\quad
j=3\cdot2^a
\]
for some integer $a\geq0$. This proves the equivalence of conditions~(1) and~(2).

The binary expansions 
$2^a=(10^a)_2$ and $3\cdot2^a=(110^a)_2$ give the addresses $L^a$ and $RL^a$, respectively, so conditions~(2) and~(3) are equivalent.  These addresses give
$\lambda(L^a)=(2^a)$ and $\lambda(RL^a)=(3\cdot2^a)$, and hence condition~(3) implies condition~(4).

Finally, condition~(4) implies condition~(1), since the
partition $(j)$ has size $j$.
\end{proof}

\begin{proposition}\label{prop:normalized-mass-gap}
If $h_j\neq j$, then
\[
\frac{h_j}{j}\leq\frac57.
\]
Equality holds if and only if $j=7\cdot2^a$ for some integer $a\geq0$.
\end{proposition}
\begin{proof} Since $\frac{h_{2j}}{2j}=\frac{h_j}{j}$, it suffices to consider odd indices. Theorem~\ref{thm:mass-fixed-points} excludes $j=1$ and $j=3$,
so write $j=2m+1\geq5$. 

If $m$ is even, then \[ h_j=h_m+1\leq m+1<\frac57(2m+1). \] If $m$ is odd, then $m\geq3$ and \[ h_j=h_m+2\leq m+2\leq\frac57(2m+1). \] 
Equality in the last inequality requires $m=3$, and $h_7=h_3+2=3+2=5$ shows that equality holds at $j=7$. Thus $7$ is the only odd index giving equality. Since doubling an index leaves the ratio unchanged, equality holds precisely when $j=7\cdot2^a$ for some integer $a\geq0$.
\end{proof}

\section{A combinatorial proof of Andrews's parity theorem}

Andrews \cite{Andrews} proved that $\sfi(n)$ is even if and only if $3\mid n$. We give a combinatorial proof by recording the weights along parent chains in the
semi-Fibonacci tree.

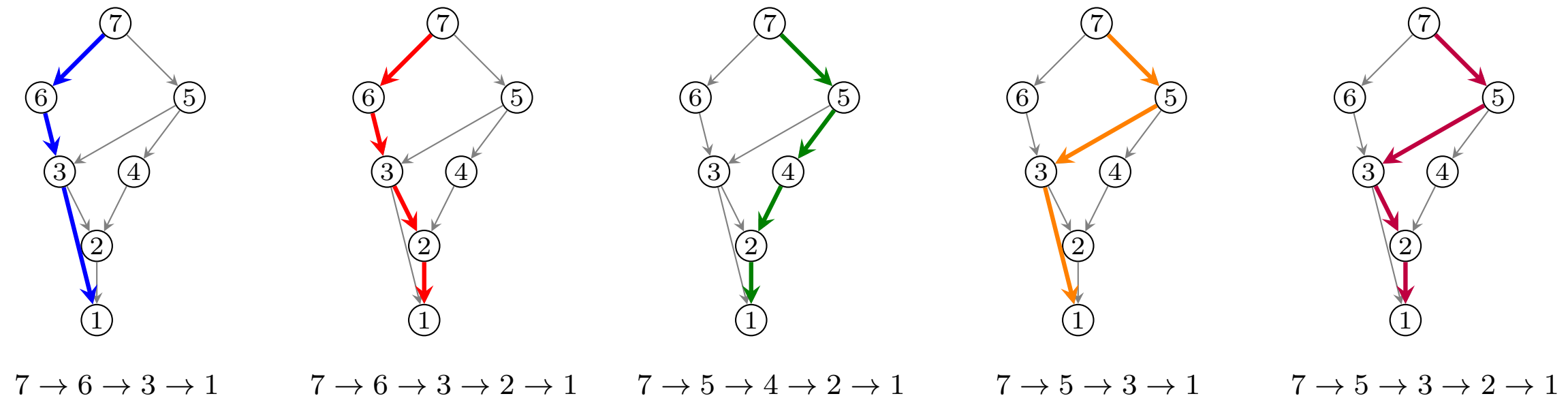
\begin{figure}[ht!]
\centering
\begin{minipage}{0.19\textwidth}
\centering
\begin{tikzpicture}[
  scale=0.9,
  >=stealth,
  every node/.style={circle, draw, inner sep=0.9pt, font=\tiny},
  base/.style={->, gray},
  hi/.style={->, blue, very thick}
]
\node (7) at (0,2.4) {$7$};
\node (6) at (-0.8,1.6) {$6$};
\node (5) at (0.8,1.6) {$5$};
\node (4) at (0.2,0.8) {$4$};
\node (3) at (-0.6,0.8) {$3$};
\node (2) at (-0.2,0) {$2$};
\node (1) at (-0.2,-0.8) {$1$};

\draw[base] (7) -- (6);
\draw[base] (7) -- (5);
\draw[base] (6) -- (3);
\draw[base] (5) -- (4);
\draw[base] (5) -- (3);
\draw[base] (4) -- (2);
\draw[base] (3) -- (2);
\draw[base] (3) -- (1);
\draw[base] (2) -- (1);

\draw[hi] (7) -- (6);
\draw[hi] (6) -- (3);
\draw[hi] (3) -- (1);
\end{tikzpicture}

\vspace{1mm}
{\scriptsize $7\to6\to3\to1$}
\end{minipage}
\hfill
\begin{minipage}{0.19\textwidth}
\centering
\begin{tikzpicture}[
  scale=0.9,
  >=stealth,
  every node/.style={circle, draw, inner sep=0.9pt, font=\tiny},
  base/.style={->, gray},
  hi/.style={->, red, very thick}
]
\node (7) at (0,2.4) {$7$};
\node (6) at (-0.8,1.6) {$6$};
\node (5) at (0.8,1.6) {$5$};
\node (4) at (0.2,0.8) {$4$};
\node (3) at (-0.6,0.8) {$3$};
\node (2) at (-0.2,0) {$2$};
\node (1) at (-0.2,-0.8) {$1$};

\draw[base] (7) -- (6);
\draw[base] (7) -- (5);
\draw[base] (6) -- (3);
\draw[base] (5) -- (4);
\draw[base] (5) -- (3);
\draw[base] (4) -- (2);
\draw[base] (3) -- (2);
\draw[base] (3) -- (1);
\draw[base] (2) -- (1);

\draw[hi] (7) -- (6);
\draw[hi] (6) -- (3);
\draw[hi] (3) -- (2);
\draw[hi] (2) -- (1);
\end{tikzpicture}

\vspace{1mm}
{\scriptsize $7\to6\to3\to2\to1$}
\end{minipage}
\hfill
\begin{minipage}{0.19\textwidth}
\centering
\begin{tikzpicture}[
  scale=0.9,
  >=stealth,
  every node/.style={circle, draw, inner sep=0.9pt, font=\tiny},
  base/.style={->, gray},
  hi/.style={->, green!50!black, very thick}
]
\node (7) at (0,2.4) {$7$};
\node (6) at (-0.8,1.6) {$6$};
\node (5) at (0.8,1.6) {$5$};
\node (4) at (0.2,0.8) {$4$};
\node (3) at (-0.6,0.8) {$3$};
\node (2) at (-0.2,0) {$2$};
\node (1) at (-0.2,-0.8) {$1$};

\draw[base] (7) -- (6);
\draw[base] (7) -- (5);
\draw[base] (6) -- (3);
\draw[base] (5) -- (4);
\draw[base] (5) -- (3);
\draw[base] (4) -- (2);
\draw[base] (3) -- (2);
\draw[base] (3) -- (1);
\draw[base] (2) -- (1);

\draw[hi] (7) -- (5);
\draw[hi] (5) -- (4);
\draw[hi] (4) -- (2);
\draw[hi] (2) -- (1);
\end{tikzpicture}

\vspace{1mm}
{\scriptsize $7\to5\to4\to2\to1$}
\end{minipage}
\hfill
\begin{minipage}{0.19\textwidth}
\centering
\begin{tikzpicture}[
  scale=0.9,
  >=stealth,
  every node/.style={circle, draw, inner sep=0.9pt, font=\tiny},
  base/.style={->, gray},
  hi/.style={->, orange, very thick}
]
\node (7) at (0,2.4) {$7$};
\node (6) at (-0.8,1.6) {$6$};
\node (5) at (0.8,1.6) {$5$};
\node (4) at (0.2,0.8) {$4$};
\node (3) at (-0.6,0.8) {$3$};
\node (2) at (-0.2,0) {$2$};
\node (1) at (-0.2,-0.8) {$1$};

\draw[base] (7) -- (6);
\draw[base] (7) -- (5);
\draw[base] (6) -- (3);
\draw[base] (5) -- (4);
\draw[base] (5) -- (3);
\draw[base] (4) -- (2);
\draw[base] (3) -- (2);
\draw[base] (3) -- (1);
\draw[base] (2) -- (1);

\draw[hi] (7) -- (5);
\draw[hi] (5) -- (3);
\draw[hi] (3) -- (1);
\end{tikzpicture}

\vspace{1mm}
{\scriptsize $7\to5\to3\to1$}
\end{minipage}
\hfill
\begin{minipage}{0.19\textwidth}
\centering
\begin{tikzpicture}[
  scale=0.9,
  >=stealth,
  every node/.style={circle, draw, inner sep=0.9pt, font=\tiny},
  base/.style={->, gray},
  hi/.style={->, purple, very thick}
]
\node (7) at (0,2.4) {$7$};
\node (6) at (-0.8,1.6) {$6$};
\node (5) at (0.8,1.6) {$5$};
\node (4) at (0.2,0.8) {$4$};
\node (3) at (-0.6,0.8) {$3$};
\node (2) at (-0.2,0) {$2$};
\node (1) at (-0.2,-0.8) {$1$};

\draw[base] (7) -- (6);
\draw[base] (7) -- (5);
\draw[base] (6) -- (3);
\draw[base] (5) -- (4);
\draw[base] (5) -- (3);
\draw[base] (4) -- (2);
\draw[base] (3) -- (2);
\draw[base] (3) -- (1);
\draw[base] (2) -- (1);

\draw[hi] (7) -- (5);
\draw[hi] (5) -- (3);
\draw[hi] (3) -- (2);
\draw[hi] (2) -- (1);
\end{tikzpicture}

\vspace{1mm}
{\scriptsize $7\to5\to3\to2\to1$}
\end{minipage}

\caption{The five directed paths in $\mathcal P_7$.}
\label{fig:P7fivecopies}
\end{figure}

For $n\geq1$, let $G_n$ be the directed graph with vertex set $\{1,2,\ldots,n\}$. For each $m\geq1$, include the edge
$2m\longrightarrow m$ if $2m\leq n$, and include the two edges $2m+1\longrightarrow2m$ and $2m+1\longrightarrow2m-1$ if $2m+1\leq n$. These edges are precisely the possible changes in weight when passing from a vertex of the semi-Fibonacci tree to its parent. Let $\mathcal P_n$ be the set of directed paths from $n$ to $1$
in $G_n$, with $\mathcal P_1$ consisting of the trivial path. A path in $\mathcal P_n$ records the successive weights along a parent chain in the semi-Fibonacci tree, starting at a partition of size $n$ and ending at the root. Figure~\ref{fig:P7fivecopies} shows the five paths in $\mathcal P_7$, each highlighted in a separate copy of $G_7$.

Recording the successive sizes along the parent chain of a partition gives a bijection from $\SF(n)$ to $\mathcal P_n$:
each path can be read backwards to reconstruct a unique partition. Hence $\sfi(n)=|\mathcal P_n|$; see
also~\cite{MorenoRamirez}.

\begin{lemma}\label{lem:canonical-path}
If $3\nmid n$, there is a unique path
$\gamma_n\in\mathcal P_n$ containing no vertex divisible by $3$.
\end{lemma}

\begin{proof}
We use induction on $n$. The case $n=1$ is immediate.  If $n$ is even,  the first step  must be  $n\longrightarrow n/2$, 
and $3\nmid n/2$. If $n>1$ is odd, the possible first steps end at $n-1$ and
$n-2$. Exactly one of these is not divisible by $3$, so there is again a unique choice. In either case, the endpoint is smaller than $n$ and not
divisible by $3$. The induction hypothesis gives a unique continuation to $1$.
\end{proof}

\begin{theorem}\label{thm:fiber-involution}
For every $n\geq1$, there is an involution
\[
\Theta_n:\mathcal P_n\longrightarrow\mathcal P_n
\]
with no fixed points when $3\mid n$ and exactly one fixed point when $3\nmid n$.
\end{theorem}

\begin{proof}
Let $P\in\mathcal P_n$. If $P$ contains no vertex divisible by $3$, set $\Theta_n(P)=P$.

Otherwise, let $m$ be the last vertex divisible by $3$ along
the path from $n$ to $1$. This vertex is odd: if it were even,
the next vertex would be $m/2$, which is also divisible by $3$,
contrary to the choice of $m$. The next step is either $m\longrightarrow m-1$ or
$m\longrightarrow m-2$. Since no vertex after $m$ is divisible
by $3$, Lemma~\ref{lem:canonical-path} shows that the path then
follows $\gamma_{m-1}$ or $\gamma_{m-2}$, respectively. Define $\Theta_n(P)$ by interchanging these two continuations
and leaving the path from $n$ to $m$ unchanged. The new path still contains no multiple of $3$ after $m$,
so $m$ remains the last vertex divisible by $3$.
Applying the operation again restores the original path.
Thus $\Theta_n$ is an involution.

Its fixed points are exactly the paths containing no vertex divisible by $3$. By Lemma~\ref{lem:canonical-path}, there is
exactly one such path, $\gamma_n$, when $3\nmid n$. When $3\mid n$, every path contains the vertex $n$, so none
is fixed.
\end{proof}

For example, when $n=7$, the involution pairs the paths as follows:
\begin{align*}
7\to6\to3\to1
&\quad\longleftrightarrow\quad
7\to6\to3\to2\to1,\\
7\to5\to3\to1
&\quad\longleftrightarrow\quad
7\to5\to3\to2\to1.
\end{align*}
The remaining path $7\to5\to4\to2\to1$ contains no multiple of $3$ and is the unique fixed point $\gamma_7$.

We next determine the difference between the numbers of partitions at even and odd depths. The length of a path is its number of edges.

\begin{theorem}\label{thm:signed-depth}
For $n\geq1$, if $P\in\mathcal P_n$ is not fixed by $\Theta_n$,
then $P$ and $\Theta_n(P)$ have lengths of opposite parity. Consequently,
\begin{equation}\label{eq:signed-depth}
\sum_{\lambda\in\mathcal{SF}(n)}
(-1)^{\operatorname{depth}(\lambda)}
=
\begin{cases}
 1, & n\equiv 1\pmod 3,\\
-1, & n\equiv 2\pmod 3,\\
 0, & n\equiv 0\pmod 3.
\end{cases}
\end{equation}
\end{theorem}

\begin{proof}
Let $j\geq1$ with $3\nmid j$. Along $\gamma_j$, the residues
modulo $3$ alternate between $1$ and $2$. Indeed, an edge
$2a\longrightarrow a$ changes the residue because $3\nmid a$,
and the endpoints of an edge from an odd vertex differ by
$1$ or $2$. Since the path ends at $1$, its length is even
when $j\equiv1\pmod3$ and odd when $j\equiv2\pmod3$.

For a nonfixed path $P$, let $m$ be its last vertex divisible
by $3$. The paths $\gamma_{m-1}$ and $\gamma_{m-2}$ have
lengths of opposite parity, since $m-1\equiv2\pmod3$ and
$m-2\equiv1\pmod3$. Each of the two tails exchanged by
$\Theta_n$ consists of one edge from $m$ followed by one
of these paths. Their lengths therefore also have opposite
parity. The path from $n$ to $m$ is unchanged, so the lengths
of $P$ and $\Theta_n(P)$ have opposite parity.

The depth of a partition is the length of its corresponding
path. Thus each nonfixed pair contributes $0$ to the signed
sum. When $3\mid n$, there are no fixed paths, and the sum
is $0$. Otherwise, only $\gamma_n$ remains. Its contribution
is $1$ when $n\equiv1\pmod3$ and $-1$ when
$n\equiv2\pmod3$, proving \eqref{eq:signed-depth}.
\end{proof}

Reducing \eqref{eq:signed-depth} modulo $2$ gives Andrews's parity
theorem.

\begin{corollary}[Andrews's parity theorem]\label{cor:parity}
For every $n>0$,
\[
\sfi(n)\equiv
\begin{cases}
0\pmod 2,&3\mid n,\\
1\pmod 2,&3\nmid n.
\end{cases}
\]
\end{corollary}

\section{The semi-Fibonacci polynomials}
The \emph{semi-Fibonacci polynomials} $S_n(x,y,z)$, introduced in our earlier work~\cite{MorenoRamirez}, are defined by
$S_0(x,y,z)=1$, $S_1(x,y,z)=y$, and, for $n\geq1$,
\[
\begin{aligned}
S_{2n}(x,y,z)&=xS_n(x,y,z),\\
S_{2n+1}(x,y,z)&=yS_{2n}(x,y,z)+zS_{2n-1}(x,y,z).
\end{aligned}
\]
In particular, $S_n(1,1,1)=\sfi(n)$. 

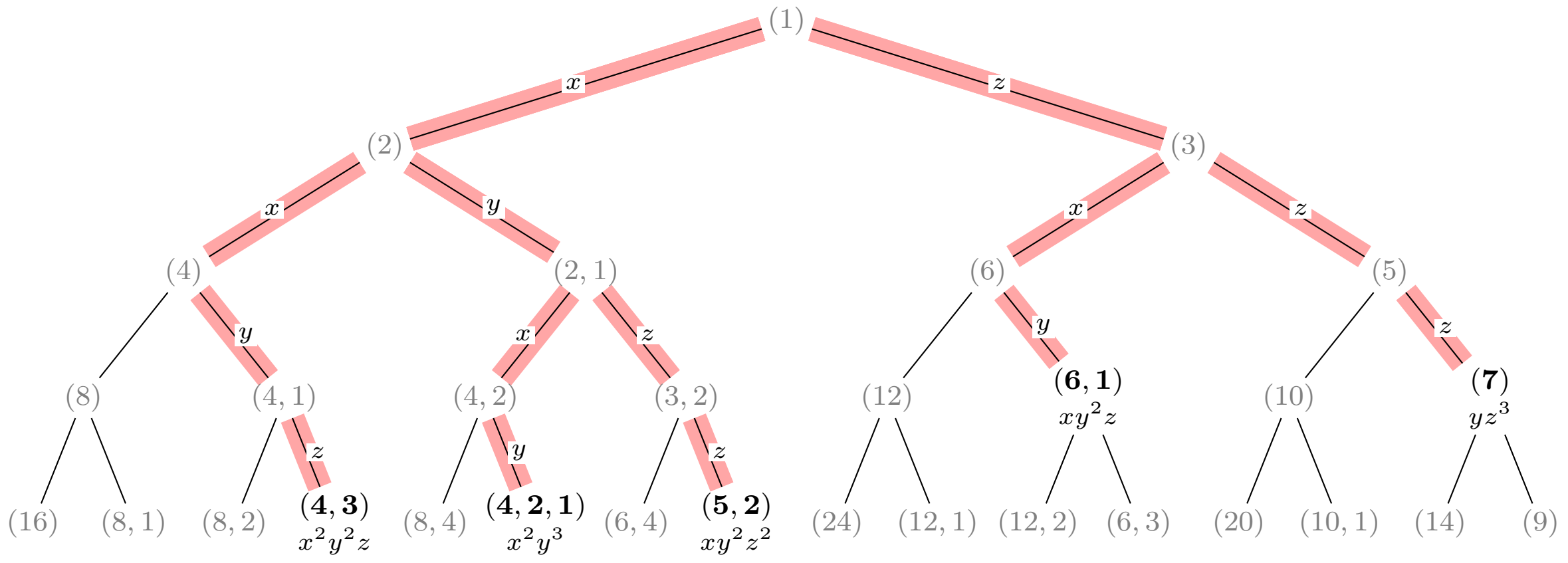
\begin{figure}[ht!]
\centering
\begin{tikzpicture}[
  level distance=1.25cm,
  every node/.style={font=\scriptsize, align=center, inner sep=1.5pt},
  edge from parent/.style={draw},
  edge label/.style={draw=none, fill=white, inner sep=1pt, font=\tiny},
  level 1/.style={sibling distance=8cm},
  level 2/.style={sibling distance=4cm},
  level 3/.style={sibling distance=2cm},
  level 4/.style={sibling distance=1cm}
]

\node (root) [text=gray] {$(1)$}
  child {node (L) [text=gray] {$(2)$}
    child {node (LL) [text=gray] {$(4)$}
      child {node [text=gray] {$(8)$}
        child {node [text=gray] {$(16)$}}
        child {node [text=gray] {$(8,1)$}}
      }
      child {node (LLR) [text=gray] {$(4,1)$}
        child {node [text=gray] {$(8,2)$}}
        child {node (LLRR) {$\mathbf{(4,3)}$\\ {\tiny $x^2y^2z$}}
          edge from parent node[edge label,right] {$z$}}
        edge from parent node[edge label,right] {$y$}}
      edge from parent node[edge label,left] {$x$}}
    child {node (LR) [text=gray] {$(2,1)$}
      child {node (LRL) [text=gray] {$(4,2)$}
        child {node [text=gray] {$(8,4)$}}
        child {node (LRLR) {$\mathbf{(4,2,1)}$\\ {\tiny $x^2y^3$}}
          edge from parent node[edge label,right] {$y$}}
        edge from parent node[edge label,left] {$x$}}
      child {node (LRR) [text=gray] {$(3,2)$}
        child {node [text=gray] {$(6,4)$}}
        child {node (LRRR) {$\mathbf{(5,2)}$\\ {\tiny $xy^2z^2$}}
          edge from parent node[edge label,right] {$z$}}
        edge from parent node[edge label,right] {$z$}}
      edge from parent node[edge label,right] {$y$}}
    edge from parent node[edge label,left] {$x$}}
  child {node (R) [text=gray] {$(3)$}
    child {node (RL) [text=gray] {$(6)$}
      child {node [text=gray] {$(12)$}
        child {node [text=gray] {$(24)$}}
        child {node [text=gray] {$(12,1)$}}
      }
      child {node (RLR) {$\mathbf{(6,1)}$\\ {\tiny $xy^2z$}}
        child {node [text=gray] {$(12,2)$}}
        child {node [text=gray] {$(6,3)$}}
        edge from parent node[edge label,right] {$y$}}
      edge from parent node[edge label,left] {$x$}}
    child {node (RR) [text=gray] {$(5)$}
      child {node [text=gray] {$(10)$}
        child {node [text=gray] {$(20)$}}
        child {node [text=gray] {$(10,1)$}}
      }
      child {node (RRR) {$\mathbf{(7)}$\\ {\tiny $yz^3$}}
        child {node [text=gray] {$(14)$}}
        child {node [text=gray] {$(9)$}}
        edge from parent node[edge label,right] {$z$}}
      edge from parent node[edge label,right] {$z$}}
    edge from parent node[edge label,right] {$z$}};

\begin{scope}[on background layer]
\draw[red!35, line width=7pt]
       (root) -- (L) -- (LL) -- (LLR) -- (LLRR);

  \draw[red!35, line width=7pt]
       (root) -- (L) -- (LR) -- (LRL) -- (LRLR);

\draw[red!35, line width=7pt]
       (root) -- (L) -- (LR) -- (LRR) -- (LRRR);

\draw[red!35, line width=7pt]
       (root) -- (R) -- (RL) -- (RLR);

 \draw[red!35, line width=7pt]
       (root) -- (R) -- (RR) -- (RRR);
\end{scope}
\end{tikzpicture}
\caption{The five paths contributing to $S_7(x,y,z)$.}
\label{fig:SF-S7}
\end{figure}

To interpret these recurrences on the tree, assign the root
$(1)$ weight $y$ and every $L$-edge weight $x$. An $R$-edge
leaving a vertex $\lambda$ has weight $y$ if $|\lambda|$ is
even and weight $z$ if $|\lambda|$ is odd.  The weight of a vertex is the product of the root weight and the edge
weights along its path from the root.  An $R$-edge has weight $y$ precisely when it immediately follows an $L$-edge, and weight $z$ otherwise. For $\lambda=\lambda(w)$, the  weight is therefore
\[
x^{|w|_L}y^{1+|w|_{LR}}z^{|w|_R-|w|_{LR}}.
\]
By Corollary~\ref{cor:word-statistics}, this is  $x^{\nu_{\max}(\lambda)}
y^{\ell(\lambda)}
z^{\sodd(\lambda)}$. Thus the tree gives a path interpretation of the partition formula from our earlier work~\cite{MorenoRamirez}:
\begin{equation}
\label{eq:SF-polynomial-enumerator}
 S_n(x,y,z)=\sum_{\lambda\in\SF(n)}
       x^{\nu_{\max}(\lambda)}y^{\ell(\lambda)}z^{\sodd(\lambda)}.
\end{equation}
For example, Figure~\ref{fig:SF-S7} highlights the five paths to partitions of $7$. The corresponding weights are 
$yz^3$, $xy^2z$, $xy^2z^2$, $x^2y^2z$, and $x^2y^3$. Hence
\[
S_7(x,y,z)
=
yz^3+xy^2z+xy^2z^2+x^2y^2z+x^2y^3.
\]

Setting $x=y=z=t$ assigns weight $t$ to the root and to every edge. A vertex at depth $d$ therefore contributes $t^{d+1}$,
so
\[
S_n(t,t,t)=
\sum_{\lambda\in\SF(n)}
t^{\operatorname{depth}(\lambda)+1}
\qquad(n\geq1).
\]
Thus the coefficient of $t^{d+1}$ in $S_n(t,t,t)$ counts the partitions of $n$ at depth $d$. For example, $S_7(t,t,t)=t(2t^3+3t^4)$, so two partitions
of $7$ occur at depth $3$ and three at depth $4$, as shown in Figure~\ref{fig:SF-S7}. We next study the levels at which
partitions of a fixed integer occur.

\subsection{The first and last levels of a fixed weight}

For $n\geq1$, define the \emph{level polynomial}
\[
\Phi_n(t)=\sum_{\lambda\in\SF(n)}t^{\depth(\lambda)}.
\]
The coefficient $[t^d]\Phi_n(t)$ counts the semi-Fibonacci partitions of $n$ at depth $d$, and $t\Phi_n(t)=S_n(t,t,t)$.  The tree rules give
\begin{equation}\label{eq:level-recurrence}
\Phi_1(t)=1,\qquad
\Phi_{2n}(t)=t\Phi_n(t),\qquad
\Phi_{2n+1}(t)=t^2\Phi_n(t)+t\Phi_{2n-1}(t).
\end{equation}
Indeed, a partition of weight $2n$ is obtained by one left step from a partition of weight $n$. A partition of weight $2n+1$ is obtained either by a left step followed by a right step from a partition of weight $n$, or by one right step from a partition of weight $2n-1$.

\begin{figure}[ht!]
    \centering
    \includegraphics[scale=0.8]{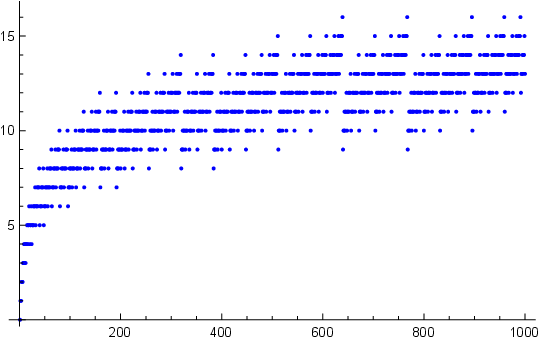}
   \qquad   \includegraphics[scale=0.8]{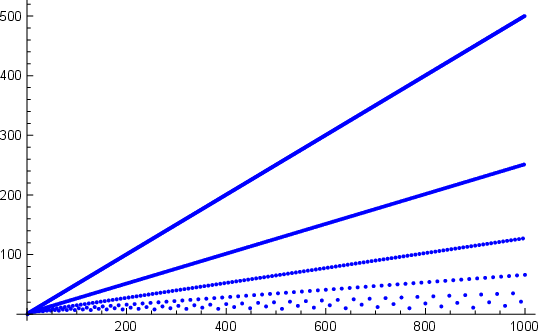}
\caption{The first and last occupied levels, $d_{\min}(n)$ and
$d_{\max}(n)$, for $1\leq n\leq 1000$.}
\label{fig:min-max-levels}
\end{figure}

Let
\[
\mathcal L(n)=\{d\geq0:[t^d]\Phi_n(t)>0\}
\]
be the set of occupied levels for weight $n$, and write
\[
d_{\min}(n)=\min\mathcal L(n),\qquad d_{\max}(n)=\max\mathcal L(n).
\]
These are the first and last levels at which weight $n$ occurs. Equivalently, they are the smallest and largest exponents with nonzero coefficients in $\Phi_n(t)$. Figure~\ref{fig:min-max-levels} plots the two sequences.

Taking the smallest exponents in \eqref{eq:level-recurrence}, we obtain
$d_{\min}(1)=0$ and, for $m\geq1$,
\begin{align*}
d_{\min}(2m)&=d_{\min}(m)+1,\\
d_{\min}(2m+1)
&=\min\{d_{\min}(m)+2,\ d_{\min}(2m-1)+1\}.
\end{align*}

\begin{lemma}\label{lem:first-level-recursion}
For $m\geq3$,
\[
d_{\min}(2m+1)=d_{\min}(m)+2.
\]
Moreover,
\[
d_{\min}(n)-d_{\min}(n-1)\leq1
\qquad(n\geq2).
\]
\end{lemma}

\begin{proof}
We prove both statements by strong induction on $n$. The initial values
\[
\begin{array}{c|rrrrrrr}
n&1&2&3&4&5&6&7\\ \hline
d_{\min}(n)&0&1&1&2&2&2&3
\end{array}
\]
verify the assertions for $n\leq7$.

First let $n=2m\geq8$. The induction hypothesis and the even recurrence give
\[
\begin{aligned}
d_{\min}(2m)-d_{\min}(2m-1)
&=d_{\min}(m)+1-\bigl(d_{\min}(m-1)+2\bigr)\\
&=d_{\min}(m)-d_{\min}(m-1)-1 \leq0.
\end{aligned}
\]
Now let $n=2m+1\geq9$. By the induction hypothesis,
\[
d_{\min}(2m-1)+1
=d_{\min}(m-1)+3
\geq d_{\min}(m)+2.
\]
The odd recurrence therefore gives
$d_{\min}(2m+1)=d_{\min}(m)+2$. Since $d_{\min}(2m)=d_{\min}(m)+1$, we also have
\[
d_{\min}(2m+1)-d_{\min}(2m)=1.
\]
This completes the induction.
\end{proof}

For a binary word $u$, let $s_2(u)$ denote its number of ones, with $s_2(\varepsilon)=0$ for the empty word. For a nonnegative
integer $n$, we also write $s_2(n)$ for its binary digit sum.

\begin{theorem}\label{thm:first-level-binary}
Let $n\geq4$ and set $k=\lfloor\log_2 n\rfloor$. Write $(n)_2=1ab\,u$, where $a,b\in\{0,1\}$ and $u$ is a binary
word of length $k-2$. Then
\[
d_{\min}(n)=k+s_2(u)+ab,
\]
where $ab$ is the ordinary product of $a$ and $b$.
\end{theorem}

\begin{proof}
The prefix $1ab$ represents the integer $p=4+2a+b$. Since
\[
d_{\min}(4)=d_{\min}(5)=d_{\min}(6)=2,
\qquad
d_{\min}(7)=3,
\]
we have $d_{\min}(p)=2+ab$.

Starting from $p$, we append the digits of $u$ from left to right. Appending a $0$ replaces $m$ by $2m$, while appending
a $1$ replaces it by $2m+1$. The even recurrence and Lemma~\ref{lem:first-level-recursion} therefore show that
each $0$ increases $d_{\min}$ by $1$ and each $1$ increases it by $2$. Since $u$ has length $k-2$, we obtain
\[
d_{\min}(n)
=2+ab+\bigl(k-2-s_2(u)\bigr)+2s_2(u)
=k+s_2(u)+ab. \qedhere
\]
\end{proof}

For example, $53=(110101)_2$, so $d_{\min}(53)=5+s_2(101)=7$. 
Theorem~\ref{thm:first-level-binary} also explains the self-similar pattern in the left panel of Figure~\ref{fig:min-max-levels}.
For $k\geq2$, the four quarters of the dyadic interval $[2^k,2^{k+1})$ correspond
to the prefixes $100$, $101$, $110$, and $111$.
The first three give $ab=0$, while the fourth gives $ab=1$. Thus the first three quarters have the same sequence of
first occupied levels, and the fourth repeats this sequence
with $1$ added to each term.

We now determine the last occupied level. Doubling all parts increases depth by $1$, so it suffices to consider odd weights.

\begin{theorem}\label{thm:max-depth}
Let $n\geq1$, and write $n=2^v u$, where $v=\nu_2(n)$ and $u=\odd(n)$. Then
\[
d_{\max}(n)=
\begin{cases}
v,&u=1,\\[1mm]
v+(u+1)/2,&u\geq3.
\end{cases}
\]
The partitions of $n$ at the last level are obtained by multiplying every part of the following partitions of $u$ by $2^v$. For $u=1$, there is only $(1)$; for $u=3$, only $(2,1)$; and for $u=5$, the two partitions $(4,1)$ and $(3,2)$.
For odd $u\geq7$, there are exactly three, with sets of parts $\{u-2,2\}$, $\{u-4,4\}$, and $\{u-6,4,2\}$.
\end{theorem}

\begin{proof}
Dividing every part by $2$ gives a bijection from $\SF(2m)$ to $\SF(m)$ and decreases depth by $1$. Repeating this $v$
times gives $d_{\max}(n)=v+d_{\max}(u)$, so it suffices to consider odd $u$.

Let $\lambda\in\SF(u)$. Write its unique odd part as $u_0$ and its even parts as $2^e u_e$, where $e$ ranges over a set
$E$ of positive integers and each $u_e$ is odd. Set $a=\max(E\cup\{0\})$. Corollary~\ref{cor:word-statistics} gives
\begin{align*}
\depth(\lambda)
&=\nu_{\max}(\lambda)+\ell(\lambda)+\sodd(\lambda)-1\\
&=a+(1+|E|)
+\left\lfloor\frac{u_0}{2}\right\rfloor
+\sum_{e\in E}\left\lfloor\frac{u_e}{2}\right\rfloor-1\\
&=a+(1+|E|)
+\frac{u_0-1}{2}
+\sum_{e\in E}\frac{u_e-1}{2}-1\\
&=a+\frac{u_0+1}{2}
+\sum_{e\in E}\frac{u_e+1}{2}-1.
\end{align*}
Since $u_0=u-\sum_{e\in E}2^e u_e$, we obtain 
\[
\depth(\lambda)
=\frac{u-1}{2}+a-
\sum_{e\in E}\frac{(2^e-1)u_e-1}{2}.
\]
If $E$ is empty, then $\lambda=(u)$ has depth  $(u-1)/2$. Otherwise, $a\geq1$. Every summand is nonnegative, and the summand corresponding to $e=a$ is at least
$\frac{(2^a-1)-1}{2}=2^{a-1}-1$. Hence
\[
\depth(\lambda)
\leq
\frac{u-1}{2}+a-(2^{a-1}-1)
\leq
\frac{u+1}{2},
\]
where the last inequality follows from $a\leq2^{a-1}$.

Equality requires $a=1$ or $a=2$. If $a=1$, equality requires $u_1=1$, so the parts are $u-2$
and $2$. For $a=2$, we have $u_2=1$, giving a part $4$.
The only other possible even part in an equality case is $2$.
Thus the sets of parts are $\{u-4,4\}$ and $\{u-6,4,2\}$.

The remaining odd part must be positive. This leaves only
$(2,1)$ when $u=3$, and $(3,2)$ and $(4,1)$ when $u=5$.
All three constructions are valid for odd $u\geq7$, and
each listed partition for $u\geq3$ has depth $(u+1)/2$.
The case $u=1$ consists only of $(1)$, of depth $0$.
This proves the result for odd weights. Multiplying every
part by $2^v$ gives the result for general $n$.
\end{proof}

For example, $72=2^3\cdot9$.  The three partitions of $9$ at the last level are $(7,2)$, $(5,4)$, and $(4,3,2)$.
Multiplying every part by $2^3$ gives the three partitions of $72$ at the last level:
\[
(56,16),\qquad (40,32),\qquad (32,24,16).
\]

The occupied levels need not form an interval. For example,
the recurrence \eqref{eq:level-recurrence} gives
\begin{align*}
\Phi_{33}(t)={}&t^6+5t^8+10t^9+8t^{10}+10t^{11}+10t^{12}\\
&+7t^{13}+6t^{14}+5t^{15}+5t^{16}+3t^{17}.
\end{align*}
Thus weight $33$ occurs at depths $6$ and $8$, but not at
depth $7$. Computing $\Phi_n(t)$ for $n<33$ shows that $33$
is the smallest weight for which $\mathcal L(n)$ is not
an interval. The next theorem shows that gaps can be
arbitrarily long.

\begin{theorem}
\label{thm:arbitrarily-long-level-gaps}
For every integer $k\geq5$,
\[
\mathcal L(2^k+1)\cap\{k+1,k+2,\ldots,2k-2\}
=\{k+1,2k-2\}.
\]
Thus weight $2^k+1$ is absent from the $k-4$ consecutive levels
\[
k+2,k+3,\ldots,2k-3.
\]
\end{theorem}
\begin{proof}
By \eqref{eq:level-recurrence},
\[
\Phi_{2^k+1}(t)
=t^2\Phi_{2^{k-1}}(t)+t\Phi_{2^k-1}(t)
=t^{k+1}+t\Phi_{2^k-1}(t).
\]
Theorem~\ref{thm:first-level-binary} gives $d_{\min}(2^k-1)=2k-3$, so the least exponent in the second summand is $2k-2$. The first summand is the single term $t^{k+1}$. Hence both degrees $k+1$ and $2k-2$ occur,
and no intermediate degree occurs.
\end{proof}

\section{Merging consecutive powers of two}
\label{sec:minimum-right-steps}

We fix the weight and minimize the number of right steps. The binary expansion describes the partitions attaining the minimum, and counting them leads to Fibonacci numbers. 

We use a bijection between semi-Fibonacci partitions and partitions into powers of $2$ in which every nonzero multiplicity is odd~\cite{Alanazi,Alegri}. Write each part uniquely as $2^i h$,
with $h$ odd, and replace it by $h$ copies of $2^i$. Distinct $2$-adic valuations ensure that every nonzero multiplicity is odd. Conversely, replacing all copies of each power of $2$ by their sum reverses
the map. For example,
\[
(6,4)=(2\cdot3,\,4\cdot1)
\quad\longleftrightarrow\quad
(4,2,2,2).
\]
In this model, an $L$-step doubles every part. An $R$-step adjoins one part $1$ if none is present, and two parts $1$ otherwise.

For $n\geq1$, write
\[
n=2^{e_1}+\cdots+2^{e_s},
\qquad e_1>\cdots>e_s\geq0.
\]
Then $\beta(n)=(2^{e_1},\ldots,2^{e_s})$ is the unique semi-Fibonacci partition of $n$ whose parts are all powers
of $2$. Starting from $\beta(n)$, choose disjoint pairs
of consecutive powers of $2$ and replace each pair by its sum:
\[
2^{e+1}+2^e\longmapsto3\cdot2^e.
\]
Here consecutive means that the exponents differ by $1$, and disjoint means that no part is used twice. Choosing
no pairs is allowed. Each merged part has the valuation of the smaller power, so the resulting parts still have distinct $2$-adic valuations and form a semi-Fibonacci partition.

For example, $\beta(15)=(8,4,2,1)$. The possible disjoint merges give
\[
\begin{array}{c|c}
\text{Pairs merged}&\text{Resulting partition}\\ \hline
\text{None}&(8,4,2,1)\\
8+4&(12,2,1)\\
4+2&(8,6,1)\\
2+1&(8,4,3)\\
8+4\text{ and }2+1&(12,3).
\end{array}
\]
The pairs $8+4$ and $2+1$ can be merged simultaneously. The pairs $8+4$ and $4+2$ cannot, since both use the part $4$.

We will show that the partitions obtained in this way are exactly those requiring the fewest right steps. For a semi-Fibonacci partition $\lambda$, let $r(\lambda)$ be
the number of right steps in its address, and define
\[
r_{\min}(n)=\min_{\lambda\in\SF(n)}r(\lambda).
\]
If $\lambda=(u_1 2^{e_1},\ldots,u_\ell2^{e_\ell})$,
with each $u_i$ odd, then
Corollary~\ref{cor:word-statistics} gives
\[
\begin{aligned}
r(\lambda)+1 =\ell(\lambda)+\sodd(\lambda)=\sum_{i=1}^{\ell}
\left(1+\left\lfloor\frac{u_i}{2}\right\rfloor\right)=\sum_{i=1}^{\ell}\frac{u_i+1}{2}.
\end{aligned}\]
In particular, merging $2^{e+1}$ and $2^e$ preserves the number of
right steps. The two original parts contribute $1+1$ to the last sum, while the merged part $3\cdot2^e$
contributes $(3+1)/2=2$.

\begin{theorem}\label{thm:minimum-right-steps}
For $n\geq1$,
\[
r_{\min}(n)=s_2(n)-1.
\]
Moreover, for $\lambda\in\SF(n)$, the following are equivalent:
\begin{enumerate}
\item $r(\lambda)=r_{\min}(n)$;
\item every part of $\lambda$ has the form $2^e$ or
      $3\cdot2^e$, and no carries occur when the parts
      are added in base $2$; 
\item $\lambda$ is obtained from $\beta(n)$ by merging disjoint
      pairs of consecutive powers of two.
\end{enumerate}
\end{theorem}

\begin{proof}
For an odd positive integer $u=2q+1$,
\[
s_2(u)=1+s_2(q)\leq1+q=\frac{u+1}{2},
\]
with equality precisely when $q\in\{0,1\}$, or equivalently $u\in\{1,3\}$. Thus, for
$\lambda=(u_12^{e_1},\ldots,u_\ell2^{e_\ell})$,
\[
r(\lambda)+1
=\sum_{i=1}^{\ell}\frac{u_i+1}{2}
\geq\sum_{i=1}^{\ell}s_2(u_i)
\geq s_2(n).
\]
The last inequality follows because each carry in binary
addition reduces the total number of ones. The partition $\beta(n)$ has $s_2(n)$ parts, with every $u_i=1$, so $r(\beta(n))+1=s_2(n)$. Hence
$r_{\min}(n)=s_2(n)-1$.

A partition attains the minimum if and only if equality holds
in both inequalities. Equality in the first means that every $u_i$ is $1$ or $3$, so every part is $2^e$ or
$3\cdot2^e$. Equality in the second means that no carries
occur when the parts are added in base $2$. This proves
the equivalence of (1) and (2).

Assume (2). Split each part $3\cdot2^e$ into $2^{e+1}$
and $2^e$. Since no carries occur, the resulting powers of $2$ are distinct. They sum to $n$, so they form $\beta(n)$. Reversing these replacements gives the
disjoint merges in (3). Conversely, merging disjoint
pairs in $\beta(n)$ produces only parts of the form
$2^e$ or $3\cdot2^e$ and introduces no carries.
Thus (2) and (3) are equivalent.
\end{proof}

For $n=15$, the five partitions listed above are therefore
exactly those with the fewest right steps. Each has
 $r(\lambda)=r_{\min}(15)=s_2(15)-1=3$.

\begin{figure}[ht!]
    \centering
    \includegraphics[scale=1]{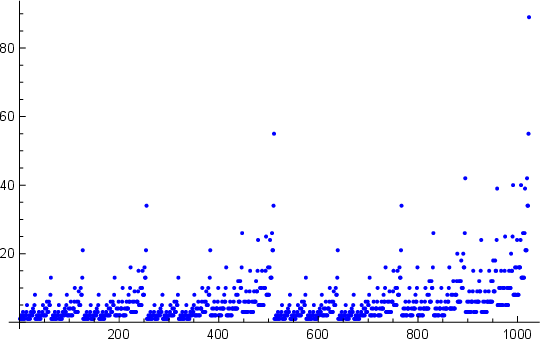}
\caption{The sequence $|\mathcal M(n)|$ for $1\leq n\leq 2^{10}-1$.}
\label{fig:M}
\end{figure}

Let
\[
\mathcal M(n)=
\{\lambda\in\SF(n):r(\lambda)=s_2(n)-1\}.
\]
Figure~\ref{fig:M} shows the first values of $|\mathcal M(n)|$.
The merging description gives the following product of
Fibonacci numbers.

\begin{proposition}\label{prop:Fibonacci-minimizers}
Write the binary expansion of $n\geq1$ as
\[
(n)_2
=
1^{a_1}0^{b_1}1^{a_2}0^{b_2}\cdots
1^{a_t}0^{b_t},
\]
where $a_i\geq1$ for $1\leq i\leq t$,
$b_i\geq1$ for $1\leq i<t$, and $b_t\geq0$. Then 
\[
|\mathcal M(n)|
=
\prod_{i=1}^{t}F_{a_i+1},
\]
where $F_{n}$ is the $n$th  Fibonacci number.
\end{proposition}

\begin{proof}
Within a block of $a_i$ ones, represent each unmerged power
of $2$ by a unit square and each merged pair by a domino. This gives a bijection with tilings of a board of length $a_i$, of which there are $F_{a_i+1}$.  Zeros prevent merges between blocks, so the choices for
different blocks are independent. Multiplying their counts
gives the formula.
\end{proof}

The product formula also explains the peaks in Figure~\ref{fig:M}. When $n=2^k-1$, its binary expansion
is a single block of $k$ ones, so $|\mathcal M(n)|=F_{k+1}$.
There is exactly one minimizing partition if and only if
the binary expansion of $n$ has no consecutive ones.
In this case no merge is possible, and the unique
minimizing partition is $\beta(n)$.

\begin{remark}
With the convention $|\mathcal M(0)|=1$, the sequence
$(|\mathcal M(n)|)_{n\geq0}$ is A246028 in the
OEIS~\cite{OEIS}. It is the run-length transform of the Fibonacci sequence \cite{WuRLT}.
Proposition~\ref{prop:Fibonacci-minimizers} therefore gives
an interpretation of A246028 in terms of partitions:
for $n\geq1$, its $n$th term counts the semi-Fibonacci
partitions of $n$ whose addresses have the fewest right steps.
\end{remark}

\section{Statistics on the semi-Fibonacci partition tree}

We now consider the joint distribution of $\nu_{\max}$, $\ell$, and $\sodd$ over the whole semi-Fibonacci partition tree.

\begin{theorem}\label{thm:global-statistics}
The joint generating function of $\nu_{\max}$, $\ell$,
and $\sodd$ is
\[
\sum_{n\geq1}S_n(x,y,z)
=
\frac{y}{(1-x)(1-z)-xy}.
\]
Equivalently, for nonnegative integers $a,b,c$, the number
of semi-Fibonacci partitions $\lambda$ satisfying
\[
\nu_{\max}(\lambda)=a,\qquad
\ell(\lambda)=b,\qquad
\sodd(\lambda)=c
\]
is
\[
\binom{a}{b-1}\binom{c+b-1}{b-1},
\]
if $1\leq b\leq a+1$, and is zero otherwise.
\end{theorem}

\begin{proof}
A part of valuation $e$ has the form $(2t+1)2^e$, with
$t\geq0$. It contributes $1$ to $\ell$ and $t$ to $\sodd$,
so its weight is $yz^t$. For a fixed valuation,
summing over all possible parts gives
\[
\sum_{t\geq0}yz^t=\frac{y}{1-z}.
\]
Now fix $\nu_{\max}(\lambda)=a$. By
Lemma~\ref{lem:valuations}, a part of valuation $a$ is
required, and each lower valuation may occur at most once.
The contribution from these partitions is therefore
\[
x^a\frac{y}{1-z}
\left(1+\frac{y}{1-z}\right)^a.
\]
Summing over $a\geq0$ gives
\[
\begin{aligned}
\sum_{n\geq1}S_n(x,y,z)
&=
\sum_{a\geq0}
x^a\frac{y}{1-z}
\left(1+\frac{y}{1-z}\right)^a\\
&=\frac{y}{(1-x)(1-z)-xy}.
\end{aligned}
\]
For the coefficient formula, suppose $1\leq b\leq a+1$.
The valuation $a$ must occur, and the remaining $b-1$
valuations can be chosen in $\binom{a}{b-1}$ ways.  Once the valuations are fixed, write the corresponding parts as
\[
(2t_1+1)2^{e_1},\ldots,(2t_b+1)2^{e_b},
\qquad t_i\geq0.
\]
The condition $\sodd(\lambda)=c$ becomes $t_1+\cdots+t_b=c$. There are
$\binom{c+b-1}{b-1}$ such weak compositions, giving the
second factor. If $b=0$ or $b>a+1$, no such partition
exists.
\end{proof}

We can also use the generating function to count partitions
by these statistics at each depth. For $d\geq0$, define
\[
A_d(x,y,z)
=
\sum_{\substack{\lambda\in V(\Tsf)\\\depth(\lambda)=d}}
x^{\nu_{\max}(\lambda)}
y^{\ell(\lambda)}
z^{\sodd(\lambda)}.
\]

\begin{corollary}\label{cor:level-enumerator}
We have
\[
\sum_{d\geq0}A_d(x,y,z)u^d
=
\frac{y}{1-(x+z)u+x(z-y)u^2}.
\]
Equivalently, $A_0=y$, $A_1=y(x+z)$, 
and, for $d\geq1$,
\[
A_{d+1}
=
(x+z)A_d+x(y-z)A_{d-1}.
\]
\end{corollary}

\begin{proof}
By Corollary~\ref{cor:word-statistics},
\[
\depth(\lambda)
=
\nu_{\max}(\lambda)+\ell(\lambda)+\sodd(\lambda)-1.
\]
Thus Theorem~\ref{thm:global-statistics} gives
\[
\begin{aligned}
\sum_{d\geq0}A_d(x,y,z)u^d
&=\frac1u\sum_{n\geq1}S_n(xu,yu,zu)\\
&=\frac{y}{1-(x+z)u+x(z-y)u^2}.
\end{aligned}
\]
The initial values and recurrence follow by multiplying
by the denominator and comparing coefficients of $u$.
\end{proof}

At depth $d$, the addresses are the $2^d$ words in
$\{L,R\}^d$. Corollary~\ref{cor:word-statistics} gives the
distributions of the individual statistics. Since
$\nu_{\max}(\lambda(w))=|w|_L$,
\[
\#\{\lambda:\depth(\lambda)=d,\ \nu_{\max}(\lambda)=a\}
=
\binom{d}{a}.
\]
Similarly, $\ell(\lambda(w))=1+|w|_{LR}$ gives
\[
\#\{\lambda:\depth(\lambda)=d,\ \ell(\lambda)=b\}
=
\binom{d+1}{2b-1}.
\]
To see this, prepend an $R$ and append an $L$ to $w$. The new word begins with $R$, ends with $L$, and has
$2b-1$ changes of letter. Choosing their positions among
the $d+1$ gaps between consecutive letters determines
the word uniquely. Finally, $\sodd(\lambda)=0$ if and only if every part of $\lambda$
is a power of $2$.  By
Corollary~\ref{cor:word-statistics}, this means that every
$R$ in the address is immediately preceded by an $L$. For $d\geq1$, these are precisely the words beginning
with $L$ and containing no consecutive $R$'s, of which
there are $F_{d+1}$. The formula also holds for $d=0$,
where the root is the only vertex. Thus exactly $F_{d+1}$ semi-Fibonacci partitions at depth $d$ have all their
parts equal to powers of $2$.

\subsection{Weights on a fixed level}

We finish by determining the average and extremal weights at each depth.

\begin{proposition}\label{prop:average-weight}
For $d\geq1$,
\[
\sum_{\substack{\lambda\in V(\Tsf)\\\depth(\lambda)=d}}
|\lambda|
=
8\cdot3^{d-1}-3\cdot2^{d-1}.
\]
Therefore, the average weight at depth $d$ is
\[
4\left(\frac32\right)^{d-1}-\frac32.
\]
\end{proposition}

\begin{proof}
Let $\Sigma_d$ be the sum of the weights at depth $d$. For $d\geq1$, half of the $2^d$ vertices have even weight and half have odd weight, since left children have even
weight and right children have odd weight.

A vertex of weight $n$ has children whose weights sum to
$3n+1$ if $n$ is even and $3n+2$ if $n$ is odd. Summing over level $d$ gives
\[
\Sigma_{d+1}=3\Sigma_d+2^d+2^{d-1},
\qquad
\Sigma_1=5.
\]
Solving the recurrence gives
\[
\Sigma_d=8\cdot3^{d-1}-3\cdot2^{d-1}.
\]
Dividing by $2^d$ gives the average.
\end{proof}

\begin{proposition}
At depth $d\geq1$, the largest weight is
$3\cdot2^{d-1}$, attained uniquely at the address $RL^{d-1}$.
For $d\geq2$, the smallest weight is $2d-1$, attained at the address $LR^{d-1}$.
\end{proposition}

\begin{proof}
The vertex with address $RL^{d-1}$ has weight $3\cdot2^{d-1}$. For a vertex of
weight $n\geq2$, the left child has weight $2n$.  The right
child has weight $n+1<2n$ if $n$ is even and $n+2<2n$ if $n$ is odd.
Since the unique maximum at depth $1$ is $3$,  induction shows that the unique maximum at each later depth is obtained by taking the left child of the preceding
maximum. This proves the first statement. For the minimum, Theorem~\ref{thm:max-depth} gives
$\depth(\lambda)\leq(|\lambda|+1)/2$, so every vertex at depth $d$ has weight at least $2d-1$.
For $d\geq2$, the partition at address $LR^{d-1}$ has parts $2$ and $2d-3$, and hence weight $2d-1$.
\end{proof}

For $d\ge2$, the smallest weight is $2d-1$, the average is
asymptotic to $\frac{8}{3}\left(\frac{3}{2}\right)^d$, and the
largest weight is $3\cdot2^{d-1}$. Thus the minimum grows linearly with depth, while the average and maximum grow exponentially, with respective
rates $3/2$ and $2$.

\section*{Declaration of AI use}
The authors used ChatGPT (OpenAI) to explore and develop
mathematical ideas they had previously formulated and to
revise the exposition. They take full responsibility for
the final manuscript.

\section*{Acknowledgments}
A. R. Moreno thanks the Universidad del Cauca for support during his Ph.D.
studies and Minciencias (Colombia) for funding through the program
``Formaci\'on de capital humano de alto nivel para promover el liderazgo
social comunitario y afrontar a la 4RI en el Departamento del Cauca,'' BPIN
2021000100129.

\end{document}